\documentclass[11pt]{article}
\usepackage[ansinew]{inputenc}
\usepackage{amssymb,hyperref}
\usepackage[dvips]{graphicx}
\hypersetup{pdfpagemode=FullScreen,colorlinks=true}

\title{$L^p$-cohomology of symmetric spaces
}
\author{Pierre Pansu$^{1,2}$\footnote{$^{1}$ Univ Paris-Sud, Laboratoire de Mathématiques d'Orsay, Orsay, F-91405 ;\hfill\eject\indent\hskip7.8pt $^{2}$ CNRS, Orsay, F-91405.}
}
\date{\today}

\newtheorem{thm}{Theorem}
\newtheorem{lem}{Lemma}
\newtheorem{prop}[lem]{Proposition}
\newtheorem{cor}[lem]{Corollary}
\newtheorem{defi}[lem]{Definition}

\def\proof{\par\medskip\noindent {\bf Proof.}{\hskip1em}}
\def\qed{~q.e.d.}

\def\R{\mathbb{R}}
\def\C{\mathbb{C}}
\def\Z{\mathbb{Z}}
\def\H{\mathbb{H}}
\def\O{\mathbb{O}}
\def\n#1{{\parallel #1 \parallel}}
\def\sp{{\rm sp}\,}
\def\dim{{\rm dim}}
\def\ker{{\rm ker}}
\def\tr{{\rm tr}}

\begin{document}
\maketitle
\begin{quote}
{\small
ABSTRACT. This is a short survey of Riemannian geometric applications of $L^p$-cohomology of thick spaces, $p\not=2$.
}
\end{quote}

\section{What is \texorpdfstring{$L^{p}$}{}-cohomology ?}
\label{intro}
\subsection{Definition}

Cohomology is the roughest invariant of topological spaces (much simpler than the fundamental group, for instance). To take a metric into account, one introduces decay conditions at infinity, this leads to $L^p$-cohomology.

\begin{defi}
\label{defcohomologie}
Let $M$ be a Riemannian manifold. Let $p>1$. Let
$L^{p}\Omega^{*}(M)$ denote the Banach space of $L^p$-differential forms and $\Omega^{*,p}(M)$ the space of $L^p$-differential forms whose exterior differential is again $L^p$, equipped with the norm
\begin{eqnarray*}
\n{\omega}_{\Omega^{*,p}}=\left(\n{\omega}_{L^p}^{p}+\n{d\omega}_{L^p}^{p}\right)^{1/p}.
\end{eqnarray*}
The complex $(\Omega^{*,p}(M),d)$ has a cohomology $H^{*,p}$ which is called the  $L^p$-{\em cohomology} of $M$. 
\end{defi}

\subsection{Reduced cohomology and torsion}

When $M$ is compact, the $L^p$ condition is no restriction, $L^p$-cohomology coincides with usual de Rham cohomology. Therefore we shall concentrate on noncompact manifolds.

In general, $H^{k,p}=(\Omega^{k,p} \cap \textrm{Ker} d)/ d\Omega ^{k-1,p}$ is not a Banach space, since $d\Omega ^{k-1,p}$ need not be closed. Therefore, one introduces a notation for its Hausdorff quotient.

\begin{defi}
\begin{eqnarray*}
R^{k,p} & = &
(\Omega^{k,p} \cap \textrm{Ker} d)/ \overline {d\Omega ^{k-1,p}},\\ 
T^{k,p}
& = & \overline {d\Omega^{k-1,p}}/ d\Omega ^{k-1,p}.
\end{eqnarray*}
$R^{k,p}$, the \emph{reduced cohomology}, is a Banach space. The
\emph{torsion} $T^{k,p}$ does not have a norm. It is not a Hausdorff
topological space.
\end{defi}

\subsection{Thick and thin ends}

For Riemannian manifolds with thin ends (e.g. Riemann surfaces with finite area cusps), one expects $L^{p}$-cohomology to coincide with the cohomology of some compactification, eventually with some small correction taking into account the specific metric behaviour of each end. This point of view has given rise to a huge litterature. We refer to S. Zucker's contribution to these proceedings.

The present notes concentrate on Riemannian manifolds with \emph{bounded geometry}, i.e. with sectional curvature bounded on both sides and injectivity radius bounded from below. For short, we call them \emph{thick spaces}. Examples of thick spaces are homogeneous spaces of Lie groups and universal coverings of compact manifolds.

We shall see that, for many thick spaces but not all, $L^{p}$-cohomology is again expressible in terms of some compactification, but the analytic (and not only topological) properties of the compactification play a role.

\subsection{Example: The real hyperbolic plane \texorpdfstring{$H^2_{\R}$}{}}

Here $H^{0,p} = 0 = H^{2,p}$ for all $p$. Let us begin with $p=2$. Since the Laplacian on $L^2$ functions is bounded below, $T^{1,2} = 0$. Therefore 
\begin{eqnarray*}
H^{1,2} &=& R^{1,2}\\ 
&=&
\{L^2\textrm{ harmonic 1-forms}\} \\
&=& \{\textrm{harmonic functions } h
\textrm{ on } H^2_{\R}\textrm{ with } \nabla h\in L^2\}/\R \,.
\end{eqnarray*}

Since the Dirichlet integral $\int \Vert \nabla h\Vert ^2$ in 2 dimensions is a
conformal invariant, one can switch from the hyperbolic metric on the disk $D$ to the euclidean metric on the disk.

Therefore
\begin{eqnarray*} 
H^{1,2} & = & \{ \textrm{harmonic functions }h\textrm{ on } D \textrm{
with } \nabla h\in L^2\}/\R\\ 
& = & \{ \textrm{Fourier series } \Sigma
a_n e^{in\theta } \textrm{ with }
  a_0 = 0, \Sigma \vert n\vert \, \vert a_n\vert ^2 < +\infty \}
\end{eqnarray*} 
which is the Sobolev space $H^{1/2}(\R/2\pi\Z)$.

More generally, for $p>1$, $T^{1,p} = 0$ and $H^{1,p}$ is equal to the Besov space
$B^{1/p}_{p,p}(\R/2\pi\Z)$ mod constants.

In this example, hyperbolic plane is compactified into a disk, $L^{p}$-cohomology identifies with a function space on the boundary circle.

\subsection{Example: the real line}

In that case, $H^{0,p}=0$. $R^{1,p}= 0$ since every function in $L^p(\R)$ can be approximated in $L^p$ with derivatives of compactly supported functions. Therefore $H^{1,p}$ is only torsion. It is non zero and 
thus infinite dimensional. Indeed, the 1-form $\frac{dt}{t}$ (cut off near the origin) is in $L^p$ for all $p>1$ but it is not the differential of a function in $L^p$.

This is an example where no compactification seems to help understanding $L^{p}$-cohomology.

\subsection{Functoriality}

By definition, $L^p$-cohomology is obviously invariant under biLipschitz diffeomorphisms. In the same way as cohomology is natural under continuous maps, and not only smooth maps, $L^p$-cohomology is natural under a wider class of maps, called \emph{uniform maps}, for thick spaces. Say a map $f:X\to Y$ between metric spaces is \emph{uniform} if $d(f(x),f(x'))$ is bounded from above in terms of $d(x,x')$ only. Proving this requires a modification of the definition, in order that it makes sense for rather general metric spaces. There are several possibilities, see \cite{Elek,Fan}. The following one is taken from \cite{Pqi}.

\begin{defi}
\label{defasympt}
Let $X$ be a metric space equipped with a measure. For each scale $t>0$, consider the simplicial complex $X_t$ whose vertices are points of $X$ and such that a subset $\Delta\subset X$ of $k+1$ points forms the vertices of a $k$-simplex if and only if its diameter is $\leq t$. Simplicial cochains $\kappa$ of $X_t$ possess an $L^p$ norm
\begin{eqnarray*}
\n{\kappa}_{p}=\left(\int_{X\times\cdots\times X}|\kappa(x_0 ,\ldots,x_k )|^p \,dx_0 \ldots dx_k \right)^{1/p}.
\end{eqnarray*}
This allows to define $L^p$-cohomology of $X_t$.
\end{defi}

\begin{prop}
\label{qiqi}
\emph{(\cite{Pqi}).} Let $X$ be a Riemannian manifold admitting a cocompact isometric group action. Assume that $H^j (X,\R)=0$ for $1\leq j\leq k$. Then, for all $t$, $H^{k,p}(X_{t})=H^{k,p}(X)$. Furthermore,
$H^{k+1,p}(X_{t})=\ker(H^{k+1,p}(X)\to H^{k+1}(X,\R))$.
\end{prop}

This means that $L^p$-cohomology depends only on the Riemannian manifold viewed as a coarse metric space. The advantage of definition \ref{defasympt} is that it applies to discrete groups as well. One merely needs to select a left-invariant metric whose balls are finite sets.

\begin{prop}
\label{qi}
Let $X$, $Y$ be measured metric spaces such that the volumes of balls are bounded above and below in terms of radius only. Then, for all $k$, for large enough $t$, there exists $t'>0$ such that $f$ induces a map $f^* : H^{k,p}(Y_{t})\to H^{k,p}(X_{t'})$. It follows that \emph{biuniform maps} induce isomorphisms on $L^p$-cohomology.
\end{prop}

Note that any homomorphism between discrete groups is uniform. In other words, $L^p$-cohomology is a functorial invariant of discrete groups which, in the case of uniform lattices in Lie groups, can probably be computed by analytic means. In this case, it is likely that the reduced part can be expressed in terms of some kind of boundary of the Lie group (compare A. Koranyi's contribution to these proceedings), but torsion cannot be excluded.

\section{What is it good for : case \texorpdfstring{$p=2$}{}
}

$L^2$ invariants constitute a rich theory by themselves, see \cite{Lueck-book}. We extract only a few illustrations of the meaning and use of $L^2$-cohomology, in connection with geometry, in the thick case. We will not develop these points, and refer to the litterature.

\subsection{Square integrable harmonic forms}

If $M$ is complete, Hodge theorem applies : $R^{k,2}$ identifies with the space of square integrable harmonic $k$-forms. Therefore, it seems to be computable by analytic means, as the example of hyperbolic plane shows. In return, $L^2$-cohomology of real, complex and quaternionic hyperbolic spaces, and especially its description in terms of differential forms on the boundary, has been used to compute $K$-theoretic invariants of the groups $SO(n,1)$, $SU(n,1)$ and $Sp(n,1)$, see \cite{Kasparov, Julg-SU, Kasparov-Julg,Julg-Sp}. 

\subsection{Torsion in degree 1}
\label{amen}

Vanishing of $T^{1,2}(M)$ is equivalent to the following isoperimetric inequality. There exists a constant $C$ such that for all smooth, compactly supported functions $u$ on $M$,
\begin{eqnarray*}
\n{u}_{L^2}\leq C\,\n{du}_{L^2}.
\end{eqnarray*}
Under local bounded geometry assumptions, this is equivalent to the following isoperimetric inequality. There exists a constant $C$ such that for all smooth, compact domains $D\subset M$,
\begin{eqnarray*}
vol(D)\leq C\,vol(\partial D).
\end{eqnarray*}
When $M$ covers a compact manifold $N$ with fundamental group $\Gamma$, isoperimetric inequality fails if and only if $\Gamma$ is \emph{amenable}, see \cite{Greenleaf}. Therefore amenable groups can be characterized by $T^{1,2}(M)\not=0$, \cite{Brooks}. 

Examples of amenable groups include $\Z$ (cf. the real line, above), solvable groups, groups of intermediate growth, see \cite{Grigorchuk-Zuk} for a state of the art.

Examples of non amenable groups include free groups, surface groups (cf. the hyperbolic plane above), nonelementary hyperbolic groups, lattices in Lie groups (except extensions of solvable Lie groups by compact Lie groups).

\subsection{Group cohomology}

Assume that $M$ is contractible and covers a compact manifold with fundamental group $\Gamma$. Then $H^{*,2}(M)$ identifies with the cohomology of the regular representation of $\Gamma$. Certain classes of groups are known to have vanising cohomology in degree 1. For instance, a finitely generated group $\Gamma$ is a \emph{Kazhdan group} if and only if $H^1 (\Gamma,\pi)=0$ for all unitary representations $\pi$ of $\Gamma$, see \cite{delaHarpe-Valette}. For such a group, $H^{1,2}(M)=0$. 

Examples of Kazhdan groups include lattices in semi-simple Lie groups with no simple factors locally isomorphic to $SO(n,1)$ or $SU(m,1)$.

Examples of non Kazhdan groups include amenable groups, free groups, surface groups and lattices in $SO(n,1)$ or $SU(m,1)$.
\begin{enumerate}
  \item Amenable groups have $T^{1,2}\not=0$ (see above) and $R^{1,2}=0$ (see \ref{bettiamen} below).
  \item Free groups and surface groups have $T^{1,2}=0$ and $R^{1,2}\not=0$.
  \item Lattices in $SO(n,1)$ or $SU(m,1)$ have $H^{1,2}=0$ except for lattices in $SO(2,1)=SU(1,1)$ which are surface groups.
\end{enumerate}

\subsection{Cohomology of towers of coverings}

Assume again that $M$ is contractible and covers a compact manifold $N$ with fundamental group $\Gamma$. If nonzero, $R^{k,2}$ is infinite dimensional. Nevertheless, one can define a kind of dimension \emph{by unit volume} $\dim_{\Gamma}(H^{k,2})$ known as the the $k$-th \emph{$L^2$-Betti number} $b^{k,2}(N)$ of $N$, \cite{Atiyah}. This works for arbitrary manifolds, but in case $\Gamma$ can be exhausted by a tower of finite index normal subgroups $\Gamma_j$, the definition is easier : renormalized usual Betti numbers $b_{k}(\Gamma_{j}\setminus M)/[\Gamma:\Gamma_{j}]$ converge to $b^{k,2}(N)$, \cite{Lueck}. Thus $R^{*,2}$ reflects the behaviour of cohomology of large compact quotients. For expositions of $L^2$-Betti numbers, see \cite{Lueck-book} and \cite{PFields}. For a connection between $L^2$-Betti numbers and spaces with thin ends, see the series of papers by J. Cheeger and M. Gromov, \cite{Cheeger-Gromov2,Cheeger-Gromov3,Cheeger-Gromov4}.

\subsection{Euler-Poincar\'e characteristic of amenable groups}
\label{bettiamen}

It has been observed early that abelian or nilpotent groups have vanishing Euler-Poincar\'e characteristic. It was not a trivial matter to extends this to solvable groups. In \cite{Cheeger-Gromov1}, J. Cheeger and M. Gromov have extended this vanishing theorem to the wide class of amenable groups. Their proof relies on $L^2$-cohomology. The additivity of $\dim_{\Gamma}$ (it adds up exactly under direct sums) implies that the Euler-Poincar\'e characteristic of a group is equal to the alternating sum of its $L^2$-Betti numbers. J. Cheeger and M. Gromov show that amenable groups have vanishing $L^2$-Betti numbers. This follows easily from the isoperimetric characterization of amenable finitely generated groups. We cheat a bit : it is a delicate point to define $L^2$-Betti numbers for arbitrary groups, not only those which admit a finite dimensional classifying space.

\subsection{Discrete series}

$R^{*,2}$ is a Hilbert space on which the isometry group of $M$ acts unitarily. In the case when $M=G/K$ is a Riemannian symmetric space, $R^{k,2}$ splits as the direct sum of irreducible representations which belong to the discrete series and have the same infinitesimal character as the trivial representation, \cite{Borel}. This provides us with a concrete realisation of these discrete series representations.

Fix a uniform lattice $\Gamma\in G$. Each a discrete series representation $\pi$ has a Harish Chandra \emph{formal dimension} $w(\pi)$, which is proportional to its $\Gamma$-dimension,
\begin{eqnarray*}
w(\pi)=\textrm{const.}(G)\frac{\dim_{\Gamma}(\pi)}{vol(\Gamma\setminus G/K)}.
\end{eqnarray*}
 $L^2$-Betti numbers are thus proportional to sums of Harish Chandra formal dimensions.

\subsection{Principal series}

Principal series representations whose infinitesimal character is equal to that of the trivial representation contribute to the torsion part of $L^2$-cohomology. Indeed, they are only weakly contained in the regular representation. This leads to the following theorem.

\begin{thm}
\label{borel}
\emph{(A. Borel, \cite{Borel}).} Let $M=G/K$ be a Riemannian symmetric space. 
\begin{enumerate}
  \item If $G$ and $K$ have equal ranks, then $T^{*,2}(M)=0$ and $R^{k,2}(M)\not=0$ if and only if $\dim (M)=2k$.
  \item Otherwise, $R^{*,2}(M)=0$ and $T^{k,2}(M)\not=0$ if and only if $2k\in(\dim(M)-\ell,\dim(M)+\ell]$, where $\ell=\mathrm{rank}_{\C}(G)-\mathrm{rank}_{\C}(K)$.
\end{enumerate}
\end{thm}

This contribution to torsion can again be quantitatively measured. Back to the general case when $M$ covers a compact manifold $N$ with fundamental group $\Gamma$. All spectral projectors $\mathbf{1}_{(0,\lambda]}(dd^{*})$ have finite $\Gamma$-traces, therefore one can define \emph{Novikov-Shubin numbers}
\begin{eqnarray*}
\alpha_{k}(N)=\limsup_{\lambda\to 0}\frac{\log\tr_{\Gamma}\mathbf{1}_{(0,\lambda]}(dd^{*})}{\log\lambda},
\end{eqnarray*}
which measure the polynomial decay of the spectral density function of $dd^{*}$ on $k$-forms, \cite{Novikov-Shubin}.

\begin{thm}
\label{olm}
\emph{(M. Olbrich, \cite{Olbrich}, N. Lohou\'e, S. Mehdi \cite{LM}).} Let $M=G/K$ be a Riemannian symmetric space. Assume $\ell=\mathrm{rank}_{\C}(G)-\mathrm{rank}_{\C}(K)>0$. Then, for all $k\in(\frac{\dim(M)-\ell}{2},\frac{\dim(M)+\ell}{2}]$, the $k$-th Novikov-Shubin invariant is equal to $\alpha_k =\ell$. Otherwise, the $L^2$-spectrum of $dd^{*}_{|\ker(d^* )^{\bot}}$ is bounded below.
\end{thm}

\section{What is it good for : case \texorpdfstring{$p\not=2$}{}
}

The theory for $p\not=2$ is much less advanced. What is missing is a generalization of $\dim_{\Gamma}$. Therefore no topological applications have been found yet. Nevertheless, there exist significant geometric applications. They are less well known, therefore we shall develop them here.

\subsection{Amenability again}

The statements of paragraph \ref{amen} extend if we replace 2 by any $p>1$.
Therefore a finitely generated group is amenable if and only if $T^{1,p}\not=0$ for some (resp. all) $p>1$.

\subsection{Euler-Poincar\'e characteristic of negatively curved manifolds}

A longstanding conjecture, attributed to H. Hopf, claims that the Euler-Poincar\'e characteristic of a compact negatively curved $2m$-manifold should be nonzero, and of the same sign as $(-1)^{m}$. J. Dodziuk and I. Singer have proposed the following approach : prove that all $L^2$-Betti numbers vanish, but the middle one, which is non zero, see \cite{Anderson,Dodziuk-Singer}. This program has been completed yet only for manifolds which admit an auxiliary K\"ahler metric, \cite{Gromovkahler}, see section \ref{hopf}.

\subsection{Hausdorff dimension at infinity}

In degree 1, $L^p$-cohomology is nondecreasing with $p$. Therefore, there is a critical $\mathbf{p}(M)$ such that $H^{1,p}(M)=0$ for $p<\mathbf{p}(M)$ and $H^{1,p}(M)\not=0$ for $p>\mathbf{p}(M)$. It turns out that in the case of hyperbolic groups, this critical exponent can be interpreted as a kind of dimension of the ideal boundary, see \cite{Pdim,Gromovasympt,Bourdon-Pajot}. 

Indeed, it is always less than or equal to the infimal Hausdorff dimension of metrics on the boundary compatible with the natural quasiconformal structure, with equality in many examples, including lattices in rank one simple Lie groups and Fuchsian hyperbolic buildings. However, the inequality may be strict, and this gives rise to examples of hyperbolic groups where the infimal Hausdorff dimension of metrics on the boundary is not achieved, see \cite{Bourdon-Pajot} and section \ref{dimension}.

\subsection{Curvature pinching}

In higher degrees, torsion in $L^p$-cohomology sometimes sharply captures negative sectional curvature pinching.

Let $-1 \leq \delta < 0$. Say a Riemannian manifold is $\delta$-pinched if its sectional curvature lies between $-a$ and $\delta a$ for some $a>0$. For example, real hyperbolic space $H^n_{\R}$ is $-1$-pinched, complex hyperbolic space $H^m_{\C}$, $m \geq 2$, quaternionic hyperbolic space $H^m_{\H}$, $m\geq 2$, and Cayley hyperbolic plane $H^2_{\O}$ are $-\frac{1}{4}$-pinched.

For real hyperbolic spaces $H^n_{\R}$, torsion in $L^p$-cohomology vanishes most of the time. In fact, in each degree, for at most one value of $p$, specificly, $p=\frac{n-1}{k-1}$ in degree $k\geq 2$. This property extends to pinched Riemannian manifolds as follows. 

\begin{thm}
\label{pinch}
\emph{\cite{Ppincement}.}
If $M^{n}$ is simply connected and $\delta$-pinched for some $\delta\in[-1,0)$, then
\begin{eqnarray*}
p<1+\frac{n-k}{k-1}\sqrt{-\delta} \quad\Rightarrow\quad T^{k,p}(M)=0.
\end{eqnarray*}
\end{thm}
This is sharp. Indeed, \emph{for every $n\geq 3$, $2\leq k\leq n-1$ and $\delta\in[-1,0)$,  there exists $\epsilon>0$ and a $\delta$-pinched homogeneous Riemannian manifold whose torsion does not vanish for $p\in (1+\frac{n-k}{k-1}\sqrt{-\delta}-\epsilon,1+\frac{n-k}{k-1}\sqrt{-\delta})$.}

However, this comparison theorem is not sharp for negatively curved symmetric spaces. For instance, for complex hyperbolic space $H^m_{\C}$, for $k=2$, the pinching comparison theorem predicts that torsion vanishes for $p<m$, whereas it turns out that torsion still does not vanish for $m\leq p<2m$, see section \ref{pincement}.

In other words, our apparently rough invariant, torsion in $L^p$-cohomology, not only detects subtle matters like sectional curvature pinching, but also distinguishes between homogeneous spaces which satisfy the same curvature bounds.

\section{Euler characteristic of negatively curved manifolds}
\label{hopf}

\subsection{Hopf's conjecture}

If $M$ is a compact Riemannian manifold with constant sectional curvature $-1$, with even dimension $n=2m$, then its Euler-Poincaré characteristic is
\begin{eqnarray*}
\chi(M)=(-1)^{m}\frac{2}{vol(S^{2m})}vol(M).
\end{eqnarray*}
Indeed, the Chern-Weil integrand $P_{\chi}(R)$ for the Euler class is a homogeneous polynomial of degree $m$ in the curvature tensor $R$. The curvature tensors of the round sphere $R_{1}$ and of real hyperbolic space $R_{-1}$ differ by a sign, $R_{-1}=-R_{1}$. Therefore $P_{\chi}(R_{-1})=(-1)^{m}P_{\chi}(R_{1})$. Integrating this relates Euler charateristics and volumes of constant curvature spaces. 

A similar calculation applies to other rank one locally symmetric manifolds $M$ (the sphere is replaced by projective spaces over the complex, quaternion and octonion numbers). In all cases, $(-1)^{m}\chi(M)>0$.

If $m=1$, the Gauss-Bonnet formula $\chi(M)=\int_{M}R$ shows that $\chi(M)<0$ as soon as $R<0$. In higher dimension, it is not so clear which negativity assumption on curvature should imply a sign for the Euler characteristic. A longstanding conjecture, attributed to H. Hopf, claims that the Euler-Poincar\'e characteristic of a compact Riemannian $2m$-manifold should be nonzero, and of the same sign as $(-1)^{m}$, if its sectional curvature is negative. 

\subsection{The Dodziuk-Singer conjecture} 

It follows readily from the definition of $L^2$-Betti numbers $b^{i,2}(M)$ that they can be used to compute Euler characteristic,
\begin{eqnarray*}
\chi(M)=\sum_{i=0}^{n}(-1)^{i}b^{i,2}(M).
\end{eqnarray*}
This has led J. Dodziuk and I. Singer to conjecture the following. Let $M$ be an $n$-dimensional compact Riemannian manifold with negative sectional curvature. Let $\tilde{M}$ denote its universal covering. Then the reduced $L^2$-cohomology $R^{i,2}(\tilde{M})$ vanishes if $i\not=n/2$, and does not vanish if $i=n/2$. 

Note that the fact that $\tilde{M}$ admits a discrete cocompact isometric group action is essential, as shown by M. Anderson, \cite{Anderson}. Indeed, there exist simply connected complete negatively curved manifolds which admit nonzero $L^{2}$-harmonic forms in several degrees simultaneously. According to \cite{Psym} Théorème G, there even exist Riemannian homogeneous spaces which admit nonzero $L^{2}$-harmonic forms in all degrees 3, 4, $\ldots,n-3$ simultaneously.

\subsection{Gromov's result}
\begin{thm}
\label{gromov}
\emph{(M. Gromov, \cite{Gromovkahler})}. Let $M$ be a compact manifold with dimension $n=2m$. Assume that $M$ admits a K\"ahler metric. Assume that the fundamental group of $M$ is Gromov-hyperbolic. Then $(-1)^{m}\chi(M)>0$.
\end{thm}

In the rest of this section, we shall give a proof of this theorem.

\subsection{Role of the K\"ahler condition}

A {\em symplectic structure} on a manifold is a differential 2-form $\omega$ such that
\begin{enumerate}
\item $\omega$ is nondegenerate, i.e. its matrix in an arbitrary basis of the tangent space is nonsingular.
\item $d\omega=0$.
\end{enumerate}
Say a Riemannian metric and a symplectic form are {\em compatible} if at each point, there exists an orthonormal basis of the tangent space in which the matrix of $\omega$ takes the block form $\left(\begin{array}{cc}0&I\\-I&0\end{array}\right)$. 

\begin{defi}
A {\em K\"ahler metric} on a manifold is the data of a Riemannian metric and parallel compatible symplectic structure.
\end{defi}

The metric and symplectic structure together determine an integrable complex structure on $M$. Clearly, complex manifolds are rare among manifolds. Therefore, admitting a K\"ahler metric is very restrictive. Nevertheless, interesting examples abound, since complex projective space admits a K\"ahler metric, which restricts to a K\"ahler metric on every smooth complex submanifold. 

The next proposition collects the two properties of K\"ahler manifolds that we shall need. See for instance \cite{Weil}, \cite{Ballmann}.

\begin{prop}
Let $(M,g,\omega)$ be a K\"ahler manifold of dimension $n=2m$. Denote by $L:\alpha\mapsto\omega\wedge\alpha$ the wedging with $\omega$ operator. Then
\begin{enumerate}
\item For all $i<m$, $L^{m-i}$ is bijection of $i$-forms to $2m-i$-forms. In particular, $L$ is injective in degrees $<m$.
\item $L$ commutes with the Laplacian on differential forms.
\end{enumerate}
\end{prop}

This implies that $L$ induces an injection on de Rham cohomology spaces $H^{i}(M,\R)\to H^{i+2}(M,\R)$ for $i<m$. Here is an $L^{2}$ variant of this statement.

\begin{cor}
\label{kahler}
Let $(M,g,\omega)$ be a complete K\"ahler manifold of dimension $n=2m$. Then wedging with $\omega$ induces an injection on reduced $L^{2}$-cohomology spaces $R^{i,2}(M)\to R^{i+2,2}(M)$ for $i<m$.
\end{cor}

\subsection{The role of negative curvature}

A simply connected nonpositively curved manifold is contractible. In fact, it possesses canonical deformation retractions : fix an origin $o$, and move an arbitrary point along the unique geodesic segment joining it to $o$. Therefore, any cycle $c$ bounds a canonical chain, the cone on $c$ with vertex $o$. Here is a variant of this construction. In a simply connected nonpositively curved manifold, from any point, one can draw a geodesic ray asymptotic to a given geodesic ray $\gamma$. This allows to construct the cone $cone_{\gamma}(c)$ on a chain $c$, as the union of all rays emanating from points of $c$, asymptotic to $\gamma$. 

\begin{center}
\includegraphics[width=3in]{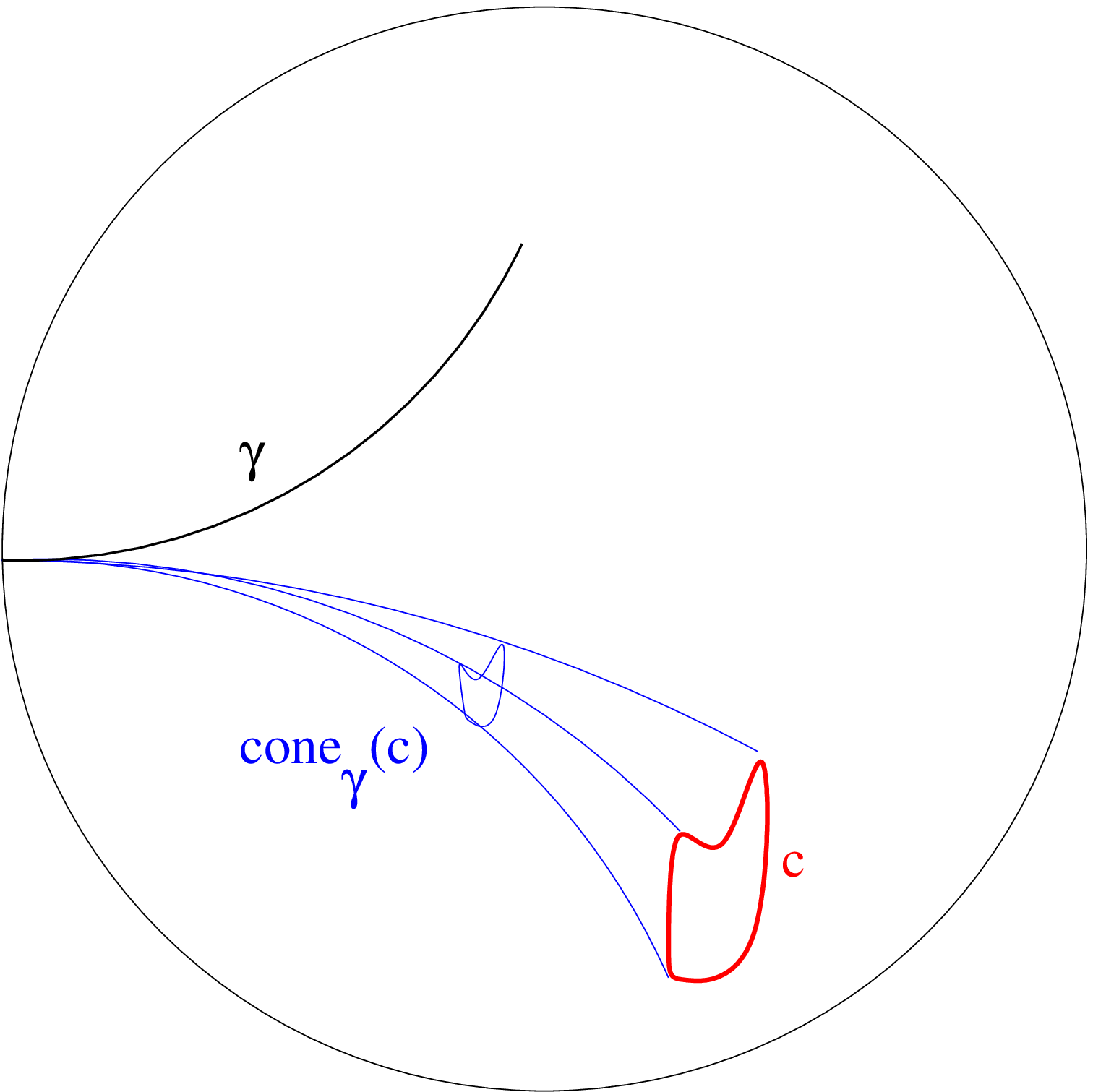}
\end{center}

We need a bound on the volume of this chain. Negative sectional curvature means that asymptotic geodesics converge exponentially. I.e., if $t\mapsto \gamma_{1}(t)$, $\gamma_{2}(t)$ are geodesics such that the distance $d(\gamma_{1}(t),\gamma_{2}(t))$ tends to zero, then
\begin{eqnarray*}
d(\gamma_{1}(t),\gamma_{2}(t))\leq e^{-\eta t}d(\gamma_{1}(0),\gamma_{2}(0)),
\end{eqnarray*}
where $-\eta^{2}$ is the upper bound of sectional curvature.

Thus, when sectional curvature is negative, the cone $cone_{\gamma}(c)$ is exponentially thin, therefore its volume is proportional to the volume of the chain,
\begin{eqnarray*}
vol(cone_{\gamma}(c))\leq \frac{1}{\eta}vol(c),
\end{eqnarray*}
provided $dim(c)\geq 1$.

Dually, the cone construction applies to differential forms. Since the norm on forms dual to the volume of chains is the $L^{\infty}$ norm, one gets

\begin{prop}
\label{cone}
{\em (M. Gromov).} Let $\tilde{M}$ be a complete simply connected Riemannian manifold with sectional curvature less than $-\eta^{2}<0$. Let $i\geq 2$.
\begin{enumerate}
\item Every $i-1$-cycle $c$ spans a $i$-chain $z$ with $vol(z)\leq \frac{1}{\eta}vol(c)$.
\item Every closed bounded differential $i$-form $\alpha$ on $\tilde{M}$ is the differential of a bounded differential $i-1$-form $\beta$ with $\n{\beta}_{L^{\infty}}\leq \frac{1}{\eta}\n{\alpha}_{L{\infty}}$.
\end{enumerate}
\end{prop}

\begin{cor}
\label{vanishk}
Assume $\tilde{M}$ covers both a compact $2m$-dimensional K\"ahler manifold and a compact negatively curved Riemannian manifold. Then
\begin{enumerate}
\item $R^{i,2}(\tilde{M})=0$ for all $i\not=m$.
\item $T^{i,2}(\tilde{M})=0$ for all $i$.
\end{enumerate}
\end{cor}

\proof
Let $\kappa$ be a reduced $L^{2}$-cohomology class, represented by a closed $L^{2}$ $i$-form $\alpha$. Since the K\"ahler form $\omega$ is bounded, Corollary \ref{cone} provides us with a bounded 1-form $\beta$ such that $\omega=d\beta$. Then $\omega\wedge\alpha=d(\beta\wedge\alpha)$, and $\beta\wedge\alpha\in L^{2}$, showing that the reduced cohomology class $[\omega\wedge\alpha]=[\omega]\wedge\kappa$ vanishes. If $i<m$, Corollary \ref{kahler} implies that $\kappa=0$. This shows that $R^{i,2}(\tilde{M})=0$.

Since
\begin{eqnarray*}
\n{d\alpha}^{2}=\langle d^{*}d\alpha,\alpha\rangle,
\end{eqnarray*}
vanishing of torsion is equivalent to a positive lower bound for the spectrum of the operator $d^{*}d$, orthogonally to its kernel, and follows from a positive lower bound for the spectrum of $\Delta=d^{*}d+dd^{*}$ orthogonally to its own kernel.

Let us show that the spectrum of the Laplacian on $i$-forms is bounded below for all $i\not=m$ (resp. orthogonally to harmonic forms in degree $i=m$). Let $\alpha$ be a $i$-form with $i>m$. Let
$\gamma$ be the the $2m-i$-form such that $\alpha=L^{i-m}\gamma$. Pointwise, $|\gamma|\leq \textrm{const.}|\alpha|$. Since $L^{i-m}$ commutes with the Laplacian, and its inverse is bounded,
\begin{eqnarray*}
\langle \Delta\gamma,\gamma\rangle\leq \textrm{const.}\langle \Delta\alpha,\alpha\rangle.
\end{eqnarray*}
Then $\alpha=d(\beta\wedge L^{i-m-1}\gamma)-\beta\wedge L^{i-m-1} d\gamma$,
\begin{eqnarray*}
\n{\alpha}^{2}&=&\langle\alpha,d(\beta\wedge L^{i-m-1} \gamma)\rangle-\langle\alpha,\beta\wedge L^{i-m-1} d\gamma\rangle\\
&\leq&\textrm{const.}(\n{d^{*}\alpha}\n{\gamma}+\n{\alpha}\n{d\gamma})\\
&\leq&\textrm{const.}\n{\alpha}\langle \Delta\alpha,\alpha\rangle.
\end{eqnarray*}
This shows that the spectrum of the Laplacian on $i$-forms is bounded below if $i>m$. Using the Hodge $*$-operator, one gets the bound for $i<m$. This gives estimates
\begin{eqnarray*}
\n{\alpha}\leq\textrm{const.}\n{d\alpha},\quad \n{\alpha}\leq\textrm{const.}\n{d^{*}\alpha},
\end{eqnarray*}
which show that $d$ and $d^{*}$ have a closed image. In particular, $dL^{2}\Omega^{m-1}(M)$ and $d^{*}L^{2}\Omega^{m+1}(M)$ are closed.

Since $\Delta$ commutes with $d$ and $d^{*}$, one gets a spectrum bound on exact and coexact $m$-forms as well. Using the Hodge decomposition
\begin{eqnarray*}
L^{2}\Omega^{m}(\tilde{M})=\textrm{ker}(\Delta)\oplus\overline{\textrm{im}(d)}\oplus\overline{\textrm{im}(d^{*})}
\end{eqnarray*}
(where closures are unneeded) gives the Laplace spectrum bound orthogonally to harmonic forms in degree $m$. This shows that $d$ has a closed image in degree $m$ as well. Therefore torsion vanishes.\qed

\subsection{Zero in the spectrum}

To complete the proof of Gromov's theorem, one needs show that reduced $L^{2}$-cohomology does not vanish in middle dimension. In view of Corollary \ref{vanishk}, if not, $L^{2}$-cohomology would vanish in all degrees. In other words, 0 would not belong to the spectrum of the Laplacian on forms in any degree. It turns out that this rarely happens for contractible spaces, as observed by J. Lott, \cite{Lott}. In fact, it follows from M. Gromov and B. Lawson's relative index theorem, \cite{Gromov-Lawson}, that 0 belongs to the spectrum for complete simply connected nonpositively curved manifolds with cover compact manifolds, see \cite{Lott}.

We cheated a bit, since Gromov's theorem has a weaker assumption : no negative curvature, only a Gromov-hyperbolic fundamental group. The proof requires two changes. First, the coning proposition has to be extended to hyperbolic metric spaces. Second, the relative index theorem has to be replaced by an adhoc avatar of the index theorem, which applies to manifolds with a cocompact action of a disconnected, nondiscrete Lie group, see \cite{Gromovkahler}, \cite{PFields}.

For variants and generalizations of Gromov's argument, see chapter 8 of \cite{Ballmann} and references therein.

\subsection{Conclusion}

The main ingredients in Gromov's argument are
\begin{enumerate}
\item The cup-product $H^{2,\infty}\otimes H^{i,2}\to H^{i+2,2}$.
\item The K\"ahler package, which gives full power to cup-product with the K\"ahler class.
\item A vanishing theorem for $L^{\infty}$-cohomology, which follows from the contracting character of canonical deformation retractions.
\end{enumerate}

\section{Dimension of the ideal boundary}
\label{dimension}

We have seen (Proposition \ref{cone}) that, for simply connected negatively curved manifolds, $L^{\infty}$-cohomology vanishes in all degrees $i\geq 2$. Note that $H^{1,\infty}$ is never zero for an unbounded space. In this section, we shall exploit the fact that, for simply connected negatively curved manifolds, $H^{1,p}$ is nonzero for $p$ large enough.

\subsection{Non vanishing of cohomology}

If $M$ is simply connected, $L^p$-cohomology in degree 1 is isomorphic to the space of functions $u$ on $M$ whose gradient is in $L^p$, modulo additive constants. Therefore, to produce non trivial $L^p$-cohomology classes, one merely needs functions whose gradient is in $L^p$, and which do not tend to a constant at infinity.

When $M$ is negatively curved, one can proceed as follows. In polar coordinates, the metric takes the form $dr^2 +g_r$, where $g_r$ is a family of Riemannian metrics on the sphere. If sectional curvature satisfies $1\leq K\leq\delta<0$, then 
\begin{eqnarray*}
(\frac{\sinh(r\sqrt{-\delta})}{\sqrt{-\delta}})^{2}g_0 \leq g_r \leq \sinh(r)^{2}g_0 ,
\end{eqnarray*}
where $g_0$ denotes the round metric on the unit sphere. In follows that the volume of the sphere of radius $r$ is at most const. $\sinh(r)^{n-1}$. Let $v$ be a smooth function on the sphere. Extend $v$ to a radial function $u$ on $M$ (i.e. $v$ does not depend on $r$). Then 
\begin{eqnarray*}
|\nabla u|\leq (\frac{\sinh(r\sqrt{-\delta})}{\sqrt{-\delta}})^{-1}.
\end{eqnarray*}
Multiply $u$ by a function of $r$ so as to make it vanish near the origin. Then 
\begin{eqnarray*}
\int_{M}|\nabla u|^{p}\leq\textrm{const.}\int_{1}^{+\infty}(\frac{\sinh(r\sqrt{-\delta})}{\sqrt{-\delta}})^{-p}\sinh(r)^{n-1}\,dr,
\end{eqnarray*}
i.e. $\nabla u\in L^p$ as soon as $p>\frac{n-1}{\sqrt{-\delta}}$. On the other hand, if $v$ is not constant, the radial extension is not in $L^p$, even up to an additive constant. Thus the cohomology class of $u$ is nonzero. We conclude that
\begin{eqnarray*}
p>\frac{n-1}{\sqrt{-\delta}}\quad\Rightarrow\quad H^{1,p}(M)\not=0.
\end{eqnarray*}
More generally, Bourdon and Pajot show that, for every nonelementary hyperbolic group, $R^{1,p}\not=0$ for $p$ large enough.

\subsection{Boundary value}
\label{vanish}

Conversely, following \cite{Strichartz}, one shows that every function $f$ on $M$ whose gradient is in $L^p$ behaves asymptotically like $u$, i.e. has a \emph{boundary value} $v_f$, i.e. a limit along almost every ray. Indeed, because of exponential growth of volume, the radial derivative $\frac{\partial f}{\partial r}$ is $L^p$ with an exponential weight, thus is $L^1$ (Hölder inequality). Clearly, if $f$ is $L^p$, $v_f$ vanishes almost everywhere. Furthermore, if $u_f$ is the radial extension of $v_f$, then $f-u_f$ is in $L^p$ (Hardy inequality, see \cite{GKS,KS}). We conclude that the map $f\mapsto v_f$ mod constants injects $H^{1,p}(M)$ into the space of $L^p$ functions on the sphere (mod constants). In particular, $H^{1,p}(M)$ is Hausdorff, and $T^{1,p}(M)=0$ (according to \ref{amen}, this is equivalent to a linear isoperimetric inequality originally due to S. T. Yau, \cite{Yau}).

\subsection{Vanishing of cohomology}

Since every function with gradient in $L^p$ is modelled on a radial function, to prove vanishing of cohomology, one merely needs check wether radial extensions have their gradient in $L^p$ or not. For negatively curved homogeneous spaces \cite{Heintze}, this happens if and only if $p>\mathbf{p}(M)=\tr(\alpha)/\min\Re e(\sp(\alpha))$, where $M$ is identified with a semi-direct product $N\times_{\alpha}\R$ defined by a derivation $\alpha$ whose eigenvalues have positive real parts. See \cite{P} for details. 

\subsection{Conformal dimension}

We give a very brief account of the idea of conformal dimension of a hyperbolic group. Precise definitions and statements should be sought in the primary litterature, see for example \cite{Bourdon-Pajot} and references therein.

A geodesic metric space is \emph{hyperbolic} if it looks like a tree : there exists a constant $\delta$ such that all triangles are $\delta$-thin, i.e., in a geodesic triangle, any point on one side is at distance at most $\delta$ of a point of one of the other two sides. Trees, of course, and hyperbolic spaces over $\R$, $\C$, $\H$ and $\O$ are typical examples.

Such metric spaces share many properties of hyperbolic spaces. They have a functorial compactification obtained by adding an \emph{ideal boundary}. Depending on the choice of an origin, a metric is defined on the ideal boundary. Changing the origin changes the metric to an equivalent one.

A finitely generated group, once a finite generating set is chosen, becomes a metric space. The group is said to be \emph{hyperbolic} if it is so as a metric space. If so, the boundary metric has a well defined Hausdorff dimension $Q$ and Hausdorff measure $\mu$, and the measure of a ball of radius $r$ is proportionnal to $r^Q$. One says that the metric is \emph{Ahlfors-$Q$-regular}. Under a change of generating set, the compactification does not change, but the metric on the ideal boundary does. It changes to a \emph{quasiconformally equivalent} one. Roughly speaking, this means that the two metrics have comparable balls, but without any control on their radii. In general, the Hausdorff dimension changes. In order to obtain a well defined numerical invariant of hyperbolic groups, one defines \emph{conformal dimension} $Cdim$ as the infimum of dimensions of metrics on the ideal boundary in the natural quasiconformal class (technical definitions vary, each of \cite{Pdim,Gromovasympt,Keith-Laakso,Bonk-Kleiner} uses a different definition).

\subsection{Dimension minimizing metrics}

The issue of wether this infimum is achieved or not is quite interesting. Based on work by S. Keith and T. Laakso, \cite{Keith-Laakso}, M. Bonk and B. Kleiner \cite{Bonk-Kleiner} show that Ahlfors $Q$-regular metrics which minimize dimension in their quasiconformal class are \emph{Loewner}. Roughly speaking, the Loewner property (introduced in \cite{Heinonen-Koskela}) means that there are enough rectifiable curves so that $Q$-capacities of condensers behave like in Euclidean space. Boundaries of negatively curved homogeneous spaces and Fuchsian buildings admit Loewner quasiconformal metrics.

G. Elek \cite{Elek-dim} and then M. Bourdon and H. Pajot \cite{Bourdon-Pajot} have extended to hyperbolic groups the discussion of \ref{vanish}. $L^p$ cohomology of the group becomes a function space $\mathcal{B}^p$ on the ideal boundary. If there exists a quasiconformal metric on the boundary with Hausdorff dimension equal to $Q$, then $\mathcal{B}^Q \not=0$. This implies that the critical exponent for $L^p$ cohomology satisfies $\mathbf{p}\leq Cdim$.

On an Ahlfors $Q$-regular Loewner space, all possible notions of dimensions, including the critical exponent for $L^p$ cohomology, coincide, \cite{Tyson,Bourdon-Pajot}. Therefore hyperbolic groups which admit dimension minimizing Ahlfors-regular quasiconformal metrics must satisfy $\mathbf{p}= Cdim$. M. Bourdon and H. Pajot construct hyperbolic amalgamations $\Gamma=A*_{C}B$ such that
\begin{eqnarray*}
Cdim(\Gamma)\geq \min\{Cdim(A),Cdim(B)\}>2,
\end{eqnarray*}
but 
\begin{eqnarray*}
L^2 b_1 (\Gamma)\geq L^2 b_1 (A)+L^2 b_1 (B)-L^2 b_1 (C)>0.
\end{eqnarray*}
Then $\mathbf{p}(\Gamma)\leq 2 < Cdim(\Gamma)$, ruining our hope of finding an optimal metric on the boundary of such a group.

\subsection{Conclusion}

The main point in this section is that, for simply connected negatively curved manifolds and groups, $H^{1,p}$ comes from the ideal boundary.

We interpret this fact as follows. Canonical deformation retractions, and specificly their variants where the vertex sits at infinity, can be used in reverse : points are moved away from the vertex instead of towards the vertex. Instead of a vanishing theorem (cohomology of space equals cohomology of point), one gets that cohomology of space equals cohomology of ideal boundary.

\section{Curvature pinching}
\label{pincement}

In this section, canonical deformation retractions and their reversed forms will be used systematicly to investigate $L^p$-cohomology in higher dimensions.

\subsection{Horospherical coordinates}

Horospherical coordinates are polar coordinates centered at a point at infinity. Such coordinates (normal exponential of a horosphere $H$) provide a diffeomorphism of $M$ onto $H\times\R$, in which the metric reads again $dr^2 +g_r$ with simpler estimates on the growth of $g_r$,
\begin{eqnarray*}
e^{2r\sqrt{-\delta}}g_0 \leq g_r \leq e^{2e}g_0 ,
\end{eqnarray*}
and the constructions of \ref{vanish} carry over \emph{verbatim}.

Denote by $\xi=\frac{\partial}{\partial r}$ the unit vector field tangent to the normals of $H$, and by $\phi_t$ its flow. The radial limit in horospherical coordinates of a function $f$ is simply $u_f =\lim_{t\to +\infty}f\circ\phi_t$. The conclusion of the discussion in \ref{vanish} can be formulated as follows. For all $p$, if $\omega=df$ is a closed $L^p$ 1-form, then the limit $u_{\omega}=du_{f}=\lim_{t\to +\infty}\phi_{t}^{*}\omega\in L_{-1}^{p}(H)$ exists, and provides an injection of $H^{1,p}(M)$ into the space of exact 1-forms on $H$ with coefficients in $L_{-1}^{p}$. Therefore $H^{1,p}(M)$ is Hausdorff. It is in this form that the method extends to higher degrees.

\subsection{Boundary values for differential forms : constant curvature}
\label{bound0}

Let $\omega$ be an $L^p$ $k$-form on $H^{n}_{\R}$. To prove that the limit $\lim_{t\to \pm\infty}\phi_{t}^{*}\omega$ exists, write
\begin{eqnarray*}
\phi_{t}^{*}\omega-\omega=d\int_{0}^{t}\phi_{s}^{*}\iota_{\xi}\omega\,ds.
\end{eqnarray*}
Since $\phi_t$ expands by a	 factor at most $e^{t}$, it expands $k-1$-forms by at most $e^{(k-1)t}$. Since $\phi_t$ expands volumes by at least $e^{(n-1)t}$, $L^p$-norms of $k-1$-forms are multiplied by at most $e^{(k-1-\frac{n-1}{p})t}$. 

If $p<\frac{n-1}{k-1}$, the integrand decreases exponentially, and the integral converges up to $+\infty$. This shows simultaneously that the limit $\omega(+\infty)$ exists, and that $\omega$ is exact if $\omega(+\infty)$ vanishes, and complete the proof that $H^{k,p}(M)$ is Hausdorff. $\omega(+\infty)$ can be interpreted as a boundary value. If furthermore $p\leq\frac{n-1}{k}$, the $L^p$ norm of $k$-forms is contracted by the reversed canonical deformation retraction, so that $\omega(+\infty)=0$ and $H^{k,p}(M)=0$. On the contrary, if $\frac{n-1}{k}<p<\frac{n-1}{k-1}$, every closed $L^p$ $k$-form on the boundary is the boundary value of some $L^p$-cohomology class, i.e. $H^{k,p}(M)\not=0$.

If $p>\frac{n-1}{k-1}$, the integrand increases exponentially, and the integral converges up to $-\infty$. But the $L^p$ norm of $k$-forms is expanded by the direct canonical deformation retraction, thus $\omega(-\infty)$ vanishes, and $\omega$ is exact. In other words, $H^{k,p}(M)=0$.

If $p=\frac{n-1}{k-1}$, none of the methods applies. Poincar\'e duality implies that reduced $L^p$-cohomology vanishes. It turns out that torsion is nonzero. These results are summed up in the following picture.

\begin{center}
\includegraphics[width=3in]{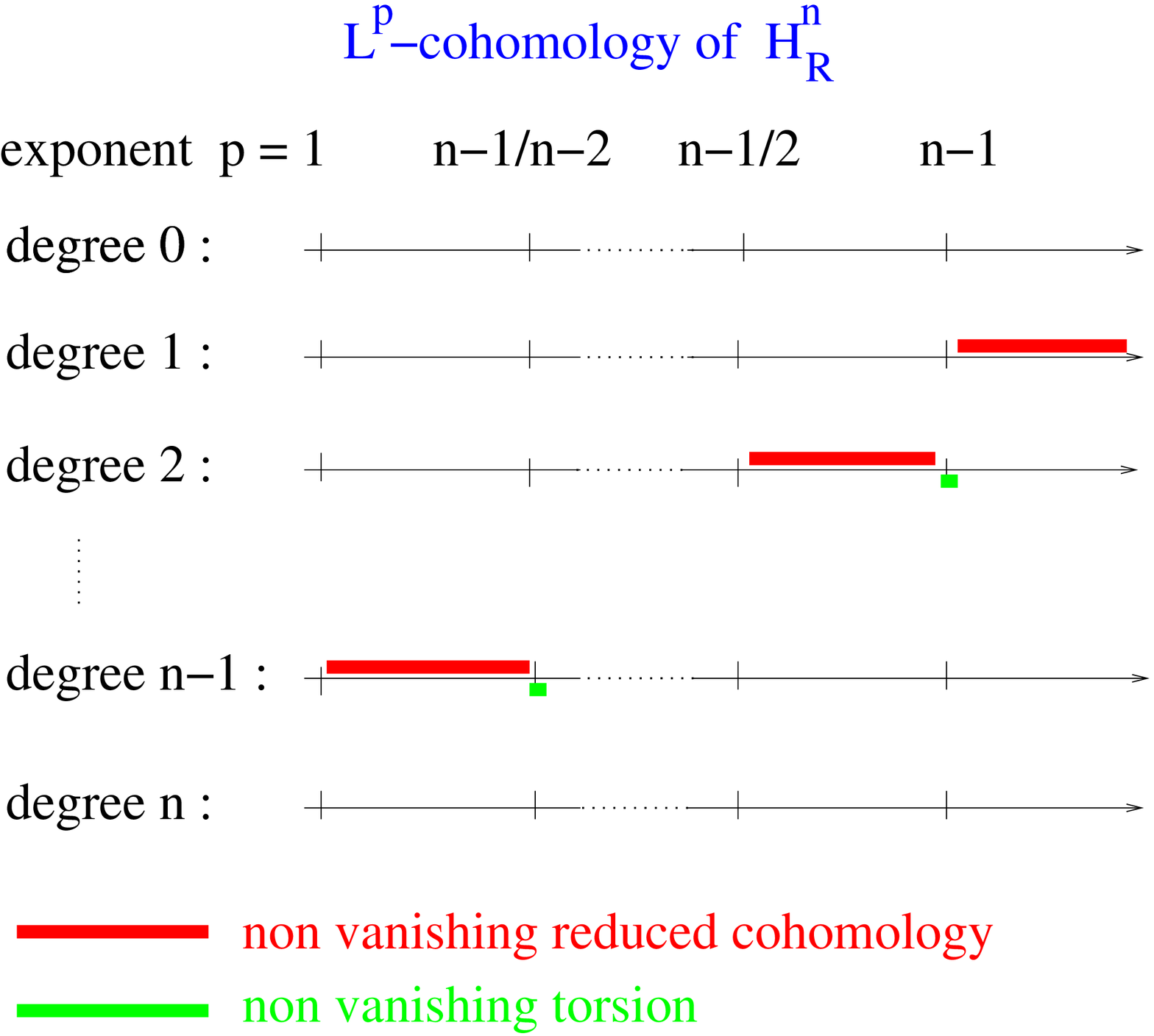}
\end{center}

\subsection{Boundary values for differential forms : pinched curvature}
\label{bound}

Assume $M$ has $\delta$-pinched curvature. Let $\omega$ be an $L^p$ $k$-form on $M$. To prove that the limit $\lim_{t\to +\infty}\phi_{t}^{*}\omega$ exists, write again
\begin{eqnarray*}
\phi_{t}^{*}\omega-\omega=d\int_{0}^{t}\phi_{s}^{*}\iota_{\xi}\omega\,ds.
\end{eqnarray*}
Again, since $\phi_t$ expands by a factor at most $e^{t}$, it expands $k-1$-forms by at most $e^{(k-1)t}$. Since $\phi_t$ expands volumes by at least $e^{(n-1)\sqrt{-\delta}t}$, $L^p$-norms of $k-1$-forms are multiplied by at most $e^{(k-1-\frac{n-1}{p}\sqrt{-\delta})t}$. If $p<\frac{n-1}{k-1}\sqrt{-\delta}$, the integrand decreases exponentially, and the integral converges up to $+\infty$. As before, this shows that $H^{k,p}(M)$ is Hausdorff.

Note that the obtained upper bound $\frac{n-1}{k-1}\sqrt{-\delta}$ is not the one stated in Theorem \ref{pinch}. A more careful inspection of the norm of the operator $\phi_{t}^{*}$ on $k-1$-form-valued $\frac{1}{p}$-densities is needed to obtain the sharp bound.

Homogeneous spaces for which Theorem \ref{pinch} is sharp are obtained as semi-direct products of $\R^{n-1}$ with $\R$ defined by diagonal matrices $\alpha$ with exactly two distinct eigenvalues. For instance, for $n=4$, $\delta=-\frac{1}{4}$ and $k=2$, the matrix $\alpha=diag(1,1,2)$ does the job. 

The same matrix defines a derivation of the Heisenberg Lie algebra $\mathrm{span}(X,Y,Z)$ with relator $[X,Y]=Z$. The corresponding semi-direct product $S=Heis\times_{\alpha}\R$ admits a left-invariant metric which makes it isometric to complex hyperbolic plane $H^{2}_{\C}$. Nevertheless, Theorem \ref{pinch} is not sharp for $2$-forms on $H^{2}_{\C}$, as we shall see next.

\subsection{Critical exponents}

Let $\xi$ be a left-invariant vectorfield on a Lie group $G$, with flow $\phi_t$. For a given degree $k\leq\dim(G)$, say an exponent $p$ is \emph{critical in degree $k$} if there exists an eigenvalue of $\Lambda^{k}ad_{\xi}$ whose real part is equal to $\tr(ad_{\xi})/p$. For instance, for $\xi=\frac{\partial}{\partial r}$ on the semi-direct product $S$ isometric to $H^{2}_{\C}$, in degree 1, there are 2 critical exponents, 2 and 4.

Below critical exponents, the flow $\phi_t$ contracts $L^p$ norms of $k$-forms, above critical exponents, $\phi_{-t}$ contracts them. In these ranges, the method of \ref{bound} applies and yields vanishing of torsion $T^{k+1,p}$. And indeed, for $H^{2}_{\C}$, Theorem \ref{pinch} predicts that $T^{2,p}=0$ for $p<2$ and for $p>4$. However, one can prove more.

\begin{thm}
\label{2}
\emph{(\cite{Psym}).}
$T^{2,p}(H^{m}_{\C})=0$ for $m\leq p<2m$.
\end{thm}

To alleviate notation, we shall treat only the 2-dimensional case.

\subsection{Two-sided boundary values}

The trick is to split each differential form $\omega$ into its contracted part $\omega_{-}$ and expanded part $\omega_{+}$, and treat each part separately. In other words, when $p$ is not critical in degree $k-1$, the operator
\begin{eqnarray*}
B:\omega\mapsto \int_{-\infty}^{0} (\phi_{t})^{*}\iota_{\xi}\omega_{-}\,dt -\int_{0}^{+\infty}(\phi_{t})^{*}\iota_{\xi}\omega_{+}\,dt,
\end{eqnarray*}
is bounded on $L^p$ $k$-forms. When $p$ is critical for no degree, $P=1-dB-Bd$ defines a chain homotopy of $\Omega^{*,p}$ to the complex $\mathcal{B}^{*,p}$ of differential forms on $G$
\begin{enumerate}
  \item whose components have $-1$ derivatives in $L^p$ ;
  \item which are killed by $\iota_{\xi}$ and $\iota_{\xi}d$.
\end{enumerate}

Although the \emph{Besov complex} $\mathcal{B}^{*,p}$ is defined in a rather implicit manner, its cohomology turns out to be computable in some cases. 

\subsection{The Besov complex for \texorpdfstring{$H^{n}_{\R}$}{}
}

The critical exponents for $H^{n}_{\R}$ are numbers of the form $\frac{n-1}{k-1}$, for $k=2,\ldots,n-1$. When $p\in(\frac{n-1}{k},\frac{n-1}{k-1})$, $\mathcal{B}^{*,p}$ vanishes in all degrees but $k$. This shows that $H^{*,p}$ is Hausdorff, is nonzero only in degree $k$ and identifies with $\mathcal{B}^{k,p}$. In that degree, it consists of differential $k$-forms whose coefficients belong to the Besov space $B_{p,p}^{-k+{{n-1}\over{p}}}(S^{n-1})$, see \cite{Triebel}.

\subsection{The Besov complex for \texorpdfstring{$H^{2}_{\C}$}{}
}

Elements of the Besov complex can be viewed as differential forms on the Heisenberg group $Heis$ with certain components missing. Denote by $(dx,dy,\tau)$ the basis of left-invariant 1-forms on $Heis$, dual to $(X,Y,Z)$. 
For $2<p<4$, a typical element of $\mathcal{B}^{1,p}$ can be written $e = f\tau $ where $f$ is a distribution on the Heisenberg group. Then
\begin{eqnarray*}
de &=& df \wedge \tau  + fd\tau \, , \\
f &=& -de(X\wedge Y)\, .
\end{eqnarray*}
One shows that there exist a constant $c$ such that
$$
\n{e}_{\mathcal{B}^p} \leq c\,\n{de}_{\mathcal{B}^p} \, .
$$
This implies that the image $d\mathcal{B}^{1,p} \subset \mathcal{B}^{2,p}
$ is closed. Therefore $T^{2,p}=0$ for $2<p<4$. This completes the proof of Theorem \ref{2}.

\subsection{Non vanishing of cohomology for \texorpdfstring{$H^{2}_{\C}$}{}}

One uses Poincar\'e duality, as formulated in \cite{GT}. Let $p'=p/p-1$ denote the conjugate exponent. Let $\omega\in\Omega^{k,p}$ be a closed form. If there exists a closed form $\psi\in\Omega^{n-k,p'}$ such that $\int_{M}\omega\wedge\psi\not=0$, then the reduced cohomology class of $\omega$ does not vanish. Otherwise, one needs to construct a sequence $\psi_{j}\in\Omega^{n-k,p'}$ such that $\n{d\psi_j}$ tends to zero but $\int_{M}\omega\wedge\psi$ does not tend to zero.

In this way, one can show that reduced cohomology is non zero in open intervals, and that torsion is non zero in two of the critical cases, $T^{2,4}$ and $T^{3,4/3}$. The remaining critical case is $k=p=2$. For this, one appeals to Theorem \ref{borel} or \cite{Gromovkahler}.

The known facts concerning the $L^p$-cohomology of $H^{2}_{\C}$ are collected in the following table. Note that although $p=2$ is a critical exponent in degree 1, nothing special seems to occur at $p=2$ on 2-forms. This is confirmed by N. Lohou\'e, \cite{Chayet-Lohoue,Lohoue}: $L^p$-cohomology varies continously in $p$ around $p=2$.

\begin{center}
\includegraphics[width=3in]{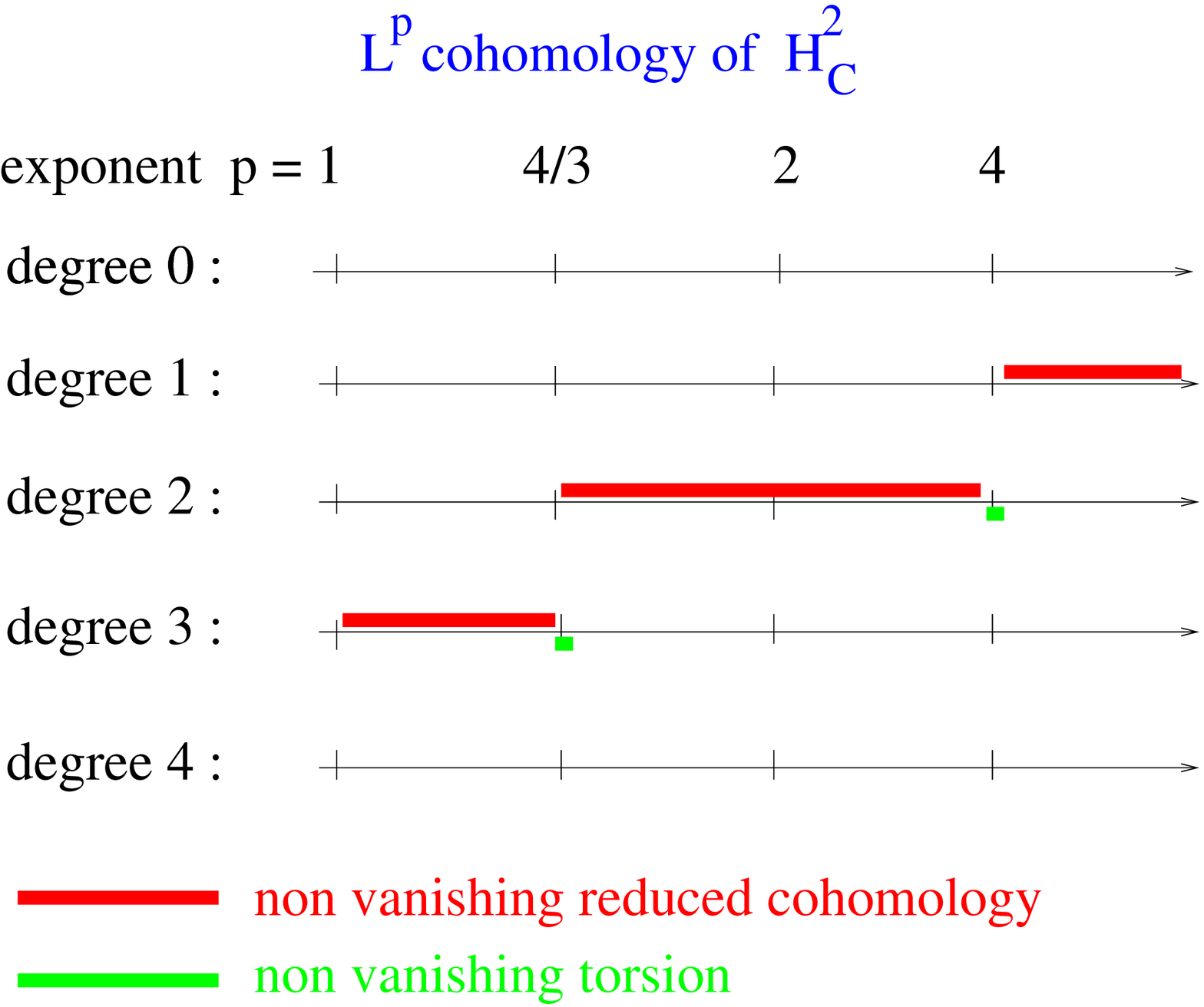}
\end{center}

The known facts concerning the $L^p$-cohomology of the $-1/4$-pinched Lie group $\R\ltimes_{diag(1,1,2)}\R^3$ are collected in the following table. It that case, the situation at $p=2$ is not fully understood.

\begin{center}
\includegraphics[width=3in]{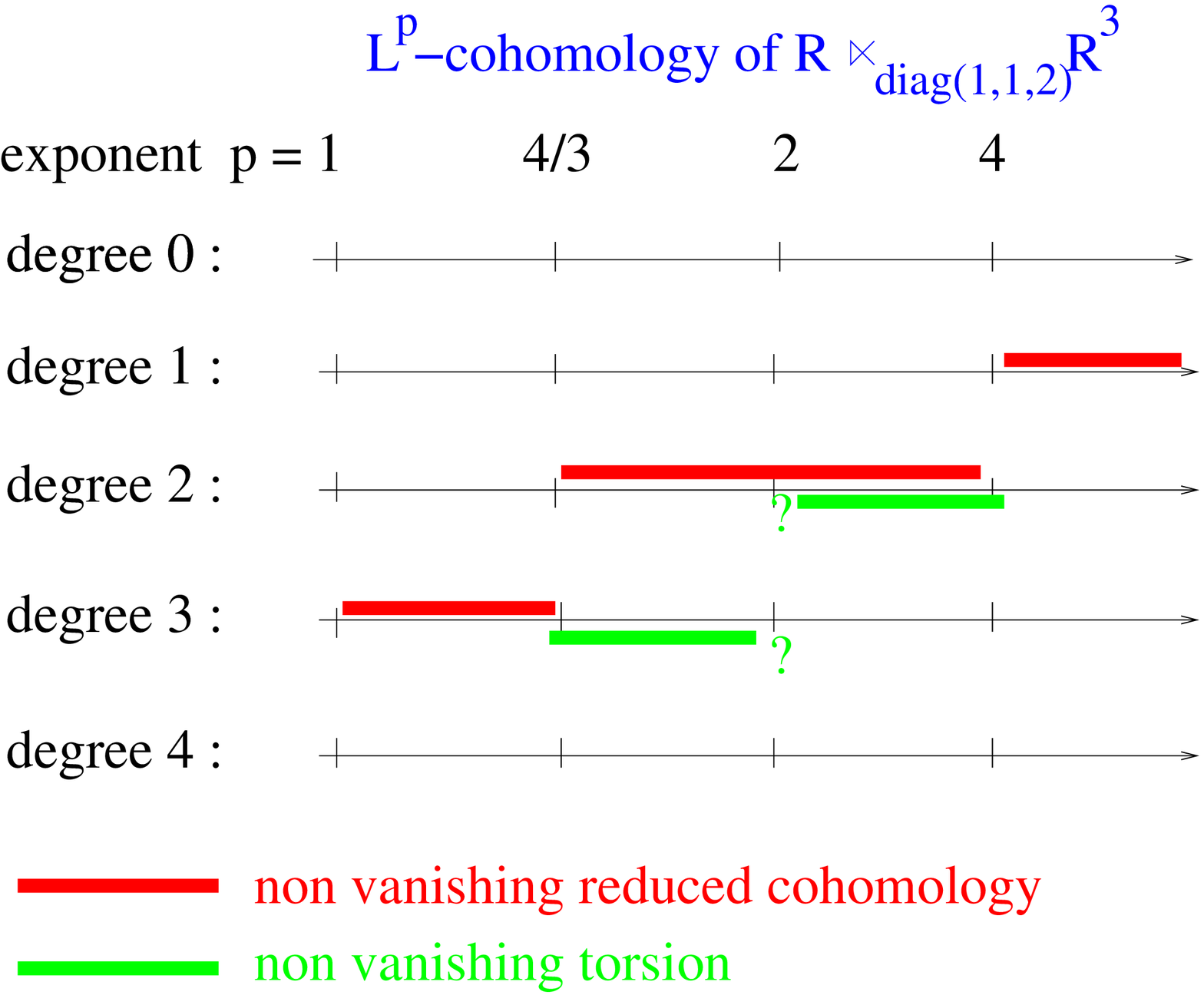}
\end{center}

\subsection{Conclusion}

All our $L^p$-cohomology calculations are based on the fact that simply connected negatively curved manifolds deformation retract onto a horosphere. The result depends on the contraction properties of this retraction, which in turn depend on the exponent $p$. For $p$ large (in degree 1, this never applies), the horosphere has to be viewed as a single point at infinity. For $p$ small (in degree 1, this applies to all $p<+\infty$), it has to be viewed as the whole ideal boundary. For intermediate, noncritical values of $p$, the retraction is used in both directions simultaneously, which is hard to interpret geometrically.

\vskip1cm
\noindent Laboratoire de Math{\'e}matique d'Orsay\\
UMR 8628 du C.N.R.S.\\
Universit{\'e} Paris-Sud XI\\
B{\^a}timent 425\\
91405 Orsay\\
France\\
\smallskip
{\tt\small Pierre.Pansu@math.u-psud.fr\\
http://www.math.u-psud.fr/$\sim$pansu}

\end{document}